\documentclass[12pt]{amsart}

\newtheorem{theorem}{Theorem}[section]

\newtheorem{corollary}[theorem]{Corollary}

\newtheorem{lemma}[theorem]{Lemma}

\theoremstyle{definition}

\newtheorem{definition}[theorem]{Definition}

\newtheorem{remark}[theorem]{Remark}

\newtheorem{example}[theorem]{Example}

\theoremstyle{parrafo}

\numberwithin{equation}{theorem}
\begin{document}

\title[]{Boundedness and unboundedness results for some maximal operators on functions of bounded variation}

\author{J. M. Aldaz and J. P\'erez L\'azaro}
\address{Departamento de Matem\'aticas y Computaci\'on,
Universidad  de La Rioja, 26004 Logro\~no, La Rioja, Spain.}
\email{aldaz@dmc.unirioja.es}
\address{Departamento de Matem\'aticas e Inform\'atica,
Universidad P\'ublica de Navarra, 31006 Pamplona, Navarra, Spain.}
\email{francisco.perez@unavarra.es}

\thanks{2000 {\em Mathematical Subject Classification.}
42B25, 26A84}

\thanks{Both authors were partially supported by Grant BFM2003-06335-C03-03 of the
D.G.I. of Spain}

\thanks{The second named author thanks the University of La Rioja for its
hospitality.}

%\subjclass{}

%\keywords{}

%\date{}

%\dedicatory{}

%\commby{}

%%% ----------------------------------------------------------------------
\begin{abstract} We  characterize the space $BV(I)$ of functions of bounded variation on an arbitrary interval  $I\subset \mathbb{R}$, in terms of a uniform boundedness condition satisfied by the  local  uncentered maximal operator $M_R$ from $BV(I)$  into the Sobolev space $W^{1,1}(I)$.
 By restriction, the corresponding characterization holds for
 $W^{1,1}(I)$. We also show that if $U$ is open in $\mathbb{R}^d, d >1$,
then boundedness from $BV(U)$  into $W^{1,1}(U)$ fails for the local
directional maximal operator $M_T^{v}$, the local strong  maximal
operator $M_T^S$, and the iterated local directional maximal
operator $M_T^{d}\circ \dots\circ M_T^{1}$. Nevertheless, if $U$
satisfies a  cone condition, then $M_T^S:BV(U)\to L^1(U)$ boundedly,
and the same happens with $M_T^{v}$, $M_T^{d} \circ \dots\circ
M_T^{1}$,  and $M_R$.
\end{abstract}

%%% ----------------------------------------------------------------------

\maketitle

%%% ----------------------------------------------------------------------

\section
{Introduction.}

\markboth{J. M. Aldaz, J. P\'erez L\'azaro}
{A characterization of $BV(I)$}

The  {\em local} uncentered Hardy-Littlewood  maximal operator $M_R$ is defined in the same way as the uncentered Hardy-Littlewood  maximal operator $M$, save for the fact that the supremum is taken over balls of diameter bounded by $R$, rather than all balls.
The terms {\em restricted} and {\em truncated} have also been used in the literature
to designate $M_R$.
We showed in \cite{AlPe}  that if $I$ is a bounded interval, then $M:BV(I)\to W^{1,1}(I)$ boundedly (Corollary 2.9). Here we complement this result by proving that for every interval
$I$, including the case of infinite length,  $M_R:BV(I)\to W^{1,1}(I)$ boundedly. Of course,
no result of this kind can hold if we consider $M$ instead of $M_R$, since $\|Mf\|_1 =\infty$
whenever $f$ is nontrivial. We shall see that if
$f\in BV(I)$,  then
$\|M_Rf\|_{W^{1,1}(I)}\le \max \{3 (1+2\log^+R), 4\} \|f\|_{BV(I)}$ (Theorem \ref{bd}), and
furthermore, the
logarithmic order of growth of
$c:= \max \{3 (1+2\log^+R), 4\}$ cannot be improved (cf. Remark \ref{log} below). Also, since $c$ is nondecreasing in $R$, it provides a uniform
bound for $M_T$ whenever $T\le R$. This observation leads to the following converse: Let
$f\ge 0$.
If there exists an $R>0$ and a constant $c = c(f,R)$ such that for all
$T\in (0,R]$,  $M_Tf\in W^{1,1}(I)$ and $\| M_Tf\|_{W^{1,1}(I)} \le c$, then $f\in BV(I)$. A fortiori, given a locally integrable $f\ge 0$,
we have that $f\in BV(I)$ if and only if
for every $R>0$, $M_Rf\in W^{1,1}(I)$ and there exists a constant $c = c(f,R)$ such that for all
$T\in (0,R]$,  $\| M_Tf\|_{W^{1,1}(I)} \le c$.  By restriction to the functions $f$ that are absolutely continuous on $I$, we obtain the corresponding characterization for $W^{1,1}(I)$.
If $f$ is real valued rather than nonnegative, since $f\in BV(I)$ (respectively $f\in W^{1,1}(I)$) if and only if both
its positive and negative parts $f^+, f^-\in BV(I)$ (respectively $f^+, f^-\in W^{1,1}(I)$),
we simply apply the previous criterion to  $M_Tf^+$ and $M_T f^-$.

 It is natural to ask whether the uniform bound condition is necessary to ensure that $f\in BV(I)$, or whether it
is sufficient just to require that  for all $T\in\mathbb{R}$,
$M_Tf\in W^{1,1}(I)$. Uniform bounds are in fact needed (see Example
\ref{counter1d}).

In higher dimensions we show that boundedness fails for the local strong maximal operator (where the supremum is taken over  rectangles with sides
parallel to the axes and uniformly bounded diameters) and
the local directional maximal operator (where the supremum is taken over  uniformly bounded segments parallel to a fixed vector), cf. Theorem \ref{Strong}  below.
But it is an open question whether the standard local  maximal operator is bounded when $d > 1$, i.e., whether given a ``sufficiently nice" open set $U\subset\Bbb R^d$, $M_R$ maps $BV(U)$ boundedly into $W^{1,1}(U)$,  or even into $BV(U)$.
On the other hand, the direction from uniform
boundedness of $M_Tf^+$ and $M_Tf^-$ to $f\in BV( U)$ follows immediately from the Lebesgue theorem on differentiation of integrals, even in the cases
of the strong and directional maximal functions (cf. Theorem \ref{trivialdir}). All the maximal operators mentioned above map $BV(U)$ boundedly into $L^1(U)$, provided
$U$ satisfies a cone condition (Theorem \ref{L1bounds}), so the question of boundedness of $M_R$ on $BV(U)$ is reduced to finding out how $DM_R$ behaves.

Previous results on these topics include the following. In \cite{Ha}, Piotr  Haj\l asz  utilized the local centered
 maximal operator to present a
 characterization, unrelated to the one given here, of the Sobolev space $W^{1,1}(\Bbb R^d)$. The boundedness of the centered Hardy-Littlewood maximal operator
on the Sobolev spaces $W^{1,p}(\Bbb R^{d})$, for $1<p\le\infty$, was proven by
Juha Kinnunen in
\cite{Ki}.  A local version of this result, valid on
$W^{1,p}(\Omega )$, $\Omega\subset\Bbb R^{d}$ open, appeared in  \cite{KiLi}. Additional
work within this line
 of research includes the papers  \cite{HaOn},  \cite{KiSa}, \cite{Lu}, \cite{Bu}, \cite{Ko1}, and
\cite{Ko2}. Of course, the case $p=1$ is significantly different from
the case $p>1$.
Nevertheless, in dimension $d=1$, Hitoshi
Tanaka showed (cf. \cite{Ta}) that if
$f\in W^{1,1}(\Bbb R)$, then the uncentered
maximal function  $Mf$ is differentiable a.e. and
$\| DMf\|_1 \le 2\| Df\|_1$ (it is asked in
\cite{HaOn}, Question 1, p. 169, whether an analogous result holds when $d > 1$).
 In \cite{AlPe} we strengthened Tanaka's result,
showing that if $f\in BV(I)$, then  $Mf$ is absolutely continuous and $\| DMf\|_1 \le | Df|(I)$, cf. \cite{AlPe}
Theorem 2.5.

  Finally we mention that the local (centered and
uncentered) maximal operator has been used in connection with inequalities involving
derivatives, cf. \cite{MaSh}  and \cite{AlPe}. Another instance of this type of
 application is given below (see Theorem \ref{ineq}).
\section{Definitions, boundedness, and unboundedness results.}

Let $I$ be an  interval  and let $\lambda$ ($\lambda^d$ if
$d > 1$) be Lebesgue measure.  Since functions of bounded variation always have lateral limits,
we can go from $(a,b)$ to $[a,b]$ by extension, and viceversa by restriction.
Thus, in what follows it does not  matter whether $I$
is open, closed or neither, nor whether it is bounded or has infinite length.
\begin{definition} We say that $f:I \to \Bbb R$ is of  bounded variation if its distributional derivative
$Df$ is a Radon measure with $|Df|(I) <\infty$, where $|Df|$ denotes the
total variation of $Df$. In higher dimensions the definition is the same, save
for the fact that $Df$ is (co)vector valued rather than real valued. More precisely,
if $U\subset \Bbb R^d$ is an open set and $f:U \to \Bbb R$ is of  bounded variation,
then  $Df$ is the vector valued Radon measure that satisfies, first,
$\int_U f \operatorname{div }\phi dx = - \int_U \phi\cdot dDf$ for all
$\phi\in C_c^1(U, \Bbb R^d)$, and  second, $|Df|(U) <\infty$.
\end{definition}

In addition to $|Df|(I) <\infty$, it is often required that $f\in L^1(I)$. We do so
only when defining the space $BV(I)$, and likewise in higher dimensions. The next definition is given only for the one dimensional case, being entirely analogous
when $d > 1$.
\begin{definition} Given the interval $I$,
$$
BV(I) := \{f:I\to \Bbb R| f\in L^1(I), Df \mbox{ is a Radon measure,
and } |Df|(I) <\infty\},
$$
and
$$
W^{1,1}(I) := \{f:I\to \Bbb R| f\in L^1(I), Df \mbox{ is a function, and } Df\in L^1(I)\}.
$$
\end{definition}

It is obvious that $W^{1,1}(I) \subset BV(I)$  properly.
The Banach space $BV(I)$ is endowed with the norm $\|f\|_{BV(I)}:= \|f\|_1 + |Df|(I)$,
and
$
W^{1,1}(I),$ with the restriction of the $BV$ norm, i.e.,  $\|f\|_{W^{1,1}(I)}:= \|f\|_1 + \|Df\|_1$.

\begin{definition} The canonical representative of $f$ is the
function
$$
\overline{f}(x) := \limsup_{\lambda (I)\to 0, x\in I}\frac{1}{\lambda (I)}\int_I f(y)dy.
$$
\end{definition}

In dimension $d=1$, bounded variation admits an elementary, equivalent definition.
Given
$P=\{x_1,\dots ,x_L\}\subset I$ with
$x_1 <\dots <x_L$, the variation of  the {\em function} $f:I\to \Bbb R$ associated to the partition $P$ is
defined as
$
V(f, I, P):= \sum_{j=2}^{L} |f(x_j) - f(x_{j-1})|,
$
and the variation of $f$ on $I$, as
$
V(f, I):=\sup_P V(f, I,P),
$
where the supremum is taken over all partitions $P$ of $I$.
Then $f$ is of bounded variation if $V(f, I) <\infty$.
As it stands this definition is not $L^p$ compatible, in the
sense that modifying $f$ on a set of measure zero can change $V(f,I)$,
and
 even make $V(f,I) = \infty.$ To remove this defect one simply says that $f$ is
of bounded variation if $V(\overline{f}, I) <\infty$. It is then well known that $|Df|(I) = V(\overline{f}, I)$.

\begin{definition}
Let $f: I\rightarrow \mathbb{R}$ be  measurable and finite a.e.. The non-increasing rearrangement $f^*$ of $f$ is defined for  $0<t<\lambda(I)$ as
\begin{equation*}
   f^*(t) = \sup_{\lambda(E)=t} \inf_{y\in E}|f(y)|.
\end{equation*}
\end{definition}

The function $f^*$ is non-increasing and equimeasurable with
$|f|$. Furthermore,
\begin{equation}\label{rear}
  \int_I f(y)dy =\int_0^{\lambda(I)}f^*(t)dt.
\end{equation}
For these and other basic properties of rearrangements see
\cite[Chapter 2]{BeSh}. We mention that the same definition can be used
for general measure spaces.

In the next definition, $\operatorname{ diam } (A)$ denotes the
diameter of a set $A$, $U\subset \Bbb R^d$ denotes an open set, and
$B\subset \Bbb R^d$ a ball with respect to some fixed norm.

\begin{definition} Given  a locally integrable function $f:U\to \Bbb R$, the {\em local}
uncentered Hardy-Littlewood
maximal function $M_R f$  is defined by
$$
M_Rf(x) := \sup_{ x\in B\subset U, \operatorname{ diam } B \le R}
\frac{1}{\lambda^d (B)}\int_B |f(y)|dy.
$$
Of course, if the bound $R$ is eliminated  then we get the usual uncentered Hardy-Littlewood
maximal function $Mf$.
\end{definition}
As noted in the introduction, the terms {\em restricted} and {\em truncated} have also been used in the literature
to designate $M_R$,
but we prefer {\em local} for the reasons detailed in Remark 2.4 of \cite{AlPe}.
Next we recall the  well known weak type (1,1) inequality satisfied by $M$ in dimension 1, with the sharp constant 2. For all $f \in L^1(I)$ and all $ t>0$,
\begin{equation}\label{weak}
  (Mf)^*(t)\le 2 \|f\|_1/t.
\end{equation}

\begin{definition}  Let  $U\subset \Bbb R^d$ be an open set, and let
 $f:U\to \Bbb R$ be a locally integrable function.
By a rectangle $R$ we mean a rectangle with sides parallel to the
axes. The local uncentered {\em strong} Hardy-Littlewood maximal
function $M_T^S f$  is defined by
$$
M_T^Sf(x) := \sup_{ x\in R\subset U, \operatorname{ diam } (R)\le T}\frac{1}{\lambda^d (R)}\int_R |f(y)|dy.
$$
Next, let $v\in \Bbb R$ be a fixed vector, and let $J$ denote a (one dimensional) segment
in $\Bbb R^d$ parallel to $v$.
The local
uncentered {\em directional} Hardy-Littlewood
maximal function $M_T^v f$  is defined by
$$
M_T^v f(x) := \sup_{ x\in J\subset U, \lambda (J)\le T}\frac{1}{\lambda (J)}\int_J |f(y)|dy.
$$
If $v=e_i$, then we write $M_T^i$ instead of $M_T^{e_i}$.
\end{definition}
We shall also be interested in the composition $M_T^d\circ\dots\circ M_T^1$ of the $d$ local directional maximal operators
in the directions of the coordinate axes, since such composition controls $M_T^S$ pointwise.
But first, we deal with the one dimensional case.
\begin{theorem}\label{bd}
    If   $|f|\in BV(I )$, then $M_Rf\in W^{1,1}(I)$ and furthermore,
$\|M_Rf\|_{W^{1,1}(I)}\le 3 (1+2\log^+ R)\|f\|_{L^1 (I)}
+4\left|D|f|\right| (I).$ Hence, $\|M_Rf\|_{W^{1,1}(I)}\le \max \{3
(1+2\log^+ R), 4\} \|f\|_{BV (I)}.$
 \end{theorem}

\begin{proof} Note that for any  interval $J$ and any $h\in BV(J)$
\begin{equation}\label{eq3}
  \|h\|_{L^\infty (J)} \le \operatorname{ essinf } |h| + |Dh| (J) \le \frac{\|h\|_{L^1(J)}}{\lambda(J)}+|Dh| (J).
\end{equation}
Now, given $f:I\to \Bbb R$, if $\left|D|f|\right|$ is a finite Radon measure on $I$, then $M_Rf$ is
absolutely continuous on $I$ and $\|DM_Rf\|_{L^1(I)}\le \left|D|f|\right|(I)$ by
\cite{AlPe}, Theorem 2.5 (we mention that for this bound on the
size of the derivative, the hypothesis $f\in L^1(I)$ is not
needed). Thus, it is enough to prove that given $|f|\in BV(I)$,
\begin{equation}\label{est}
  \|M_Rf\|_{L^1(I)}\le 3(1+2\log^+R) \|f\|_{L^1(I)} + 3 |D|f||(I).
\end{equation}
We may assume that $0\le f=\bar{f}$, since this does not change
any value of $M_Rf$. Given $k\in\mathbb{Z}$ we denote by $I_k$ and $J_k$
the (possibly empty) intervals $I\cap[kR,(k+1)R)$ and $I\cap[(k-1)R,(k+2)R)$ respectively.
We also set $f_k := f|_{J_k}$. Fix $k$. Then

\begin{equation}\label{eq2}
  \int_{I_k}M_Rf(x)dx =
  \int_{I_k}M_Rf_k(x)dx \le
  \int_{I_k}Mf_k(x)dx.
\end{equation}
Suppose first that $\lambda(I_k)\le 1$. From (\ref{eq3}) we get
\begin{equation}\label{smalllamb}
  \int_{I_k}Mf_k(x)dx \le   \lambda (I_k) \|f_k\|_{L^\infty (J_k)} \le
  \|f_k\|_{L^1(J_k)}+|Df_k| (J_k).
\end{equation}
And if $\lambda (I_k) > 1$, then from (\ref{rear}) and (\ref{weak}) we obtain
\begin{equation}\label{biglamb}
  \int_{I_k}Mf_k(x)dx = \int_0^{\lambda (I_k)} (Mf_k)^*(t)dt
  = \int_0^1 +\int_1^{\lambda (I_k)}
\end{equation}
\begin{equation*}
\le\|f_k\|_{L^\infty (J_k)} +2 \|f_k\|_{L^1 (J_k)}\int_1^{\lambda (I_k)}t^{-1}dt
\end{equation*}
\begin{equation*}
\le (1+2\log R)\|f_k\|_{L^1 (J_k)} +|Df_k| (J_k).
\end{equation*}
Since the intervals $I_k$ are all disjoint, and each nonempty $I_k$ is contained in
$J_{k-1}, J_k$ and $J_{k+1}$, having empty intersection with all  the other $J_i$'s, the
 estimates (\ref{smalllamb}) and (\ref{biglamb})  yield
\begin{equation}\label{L1norm}
\|M_R f\|_{L^1(I)} =  \sum_{-\infty}^{\infty} \int_{I_k}M_Rf(x)dx
\end{equation}
\begin{equation*}
\le \sum_{-\infty}^{\infty} \left((1+2\log^+ R)\|f_k\|_{L^1 (J_k)}
+|Df_k| (J_k)\right)  \end{equation*}
\begin{equation*}
= 3\sum_{-\infty}^{\infty} (1+2\log^+ R)\|f_k\|_{L^1 (I_k)}
+3\sum_{-\infty}^{\infty}|Df_k| (I_k)  \end{equation*}
\begin{equation*}
= 3 (1+2\log^+ R)\|f\|_{L^1 (I)} +3|Df| (I).
\end{equation*}
Thus,
\begin{equation*}\label{Rest}
\|M_R f\|_{BV(I)}  \le 3 (1+2\log^+ R)\|f\|_{L^1 (I)} +4|Df| (I) \le
\max \{3 (1+2\log^+ R), 4\} \|f\|_{BV (I)}.
\end{equation*}
\end{proof}

\begin{remark}\label{log} The example $f:\Bbb R \to \Bbb R$ given by $f:=\chi_{[0,1]}$ shows that
the logarithmic order of growth  in  the preceding theorem is the correct
one. Here all the relevant quantities can be easily computed: $\|f\|_{L^1 (\Bbb R)}= 1$,
$|Df|(\Bbb R) = 2$, $\|M_R f\|_{L^1 (\Bbb R)}= 1 + 1/ R + 2\log R$ for $R\ge 1$, and
$|DM_Rf|(\Bbb R) = 2$ (for all $R>0$).
\end{remark}

As noted in \cite{AlPe}, this kind of bounds on the size of maximal functions and their
derivatives can be
used to obtain variants of the classical Poincar\'e inequality, as well as other inequalities involving
derivatives,  under less regularity,
 by using  $DM_R f$ (a function) instead
of  $Df$ (a Radon measure).
Here we present another instance of the same idea, a  Poincar\'e type inequality involving $\|M_Rf\|_1$; the
argument is standard but short, so we include it for the reader's convenience.

 \

Given a compactly supported function $f$, denote by
$N(f,R):= \operatorname{ supp }f + [-R, R]\subset \Bbb R$ the closed $R$-neighborhood
of its support, that is, the set of all points at distance less than or equal to $R$ from
the support of $f$.

 \begin{theorem}\label{ineq}
    Let   $f\in BV(\Bbb R )$ be compactly supported. Then for all $R> 0$, we have
$ \|f\|_2^2 $
 \begin{equation*}\label{MRpoin}
      \le \min\left\{ \frac{(3 (1+2\log^+R))^2}{\lambda (N(f,R))} \|f\|_{BV(\mathbb{R})}^2+
 \left(\frac{ \left(\lambda(N(f,R))\right)^2}{2}\right) \|D M_Rf\|_2^2, \lambda(N(f,R))^2
 \|D M_Rf\|_2^2, \right\}.
      \end{equation*}
 \end{theorem}
\begin{proof} Let $x< y$ be points in $\Bbb R$. By the
Fundamental Theorem of Calculus,
\begin{equation*}
M_R f(y) - M_R f(x)  = \int_x^y DM_Rf(t) dt \le \|D M_Rf\|_1.
      \end{equation*}
Squaring and integrating with respect to $x$ and $y$ over
$N(f,R)^2$, we get
\begin{equation*}
      \|M_Rf\|_2^2  \le
\frac{\|M_Rf\|_1^2}{\lambda (N(f,R))}+ \|D M_Rf\|_1^2\left( \frac{\lambda (N(f,R))}{2}\right).
      \end{equation*}
Since $ \|f\|_2^2  \le  \|M_Rf\|_2^2$, using  (\ref{L1norm}) and either Jensen or H\"older inequality we obtain
 \begin{equation*}\label{b}
      \|f\|_2^2  \le  \frac{(3 (1+2\log^+R))^2}{\lambda (N(f,R))} \|f\|_{BV(\mathbb{R})}^2+
 \left(\frac{ \left(\lambda(N(f,R))\right)^2}{2}\right) \|D M_Rf\|_2^2.
      \end{equation*}
On the other hand, integrating $M_R f(y)   = \int_{\infty}^y DM_Rf(t) dt \le \|D M_Rf\|_1$
and repeating the previous steps we get
 \begin{equation*}\label{a}
      \|f\|_2^2  \le \lambda(N(f,R))^2
 \|D M_Rf\|_2^2.
      \end{equation*}
\end{proof}
 \begin{remark} In connection with the preceding inequality, we point out that if
$1 < p < \infty$ and $f\in W^{1,p}(\Bbb R )$, then $\|D M_Rf\|_p\le c_p \|D f\|_p$,
with $c_p$ independent of $R$. Of course, the interest of the result lies in the
fact that we can have  $\|D M_Rf\|_p < \infty$ even if $Df$ is not a function (standard
example, $f=\chi_{[0,1]}$). The cases $p=1,\infty$ are handled in \cite{AlPe}, Theorems
2.5 and 5.6. There we have $\|D M_Rf\|_p\le  \|D f\|_p$. To see why $\|D M_Rf\|_p\le c_p \|D f\|_p$ holds
with $c_p$ independent of $R$, repeat
the sublinearity argument from \cite{Ki}, Remark 2.2 (i) (cf. also
\cite{HaOn}, Theorem 1) using  $M_Rf\le Mf$ to remove the dependency of the constant on
$R$.
\end{remark}

We shall consider next  the local strong, directional, and iterated
directional maximal operators, proving boundedness from $BV(U)$ into
$L^1(U)$ and lack of boundedness from $BV(U)$ into $BV(U)$. Of
course, since the strong maximal operator dominates pointwise (up to
a constant factor) the maximal operator associated to an arbitrary
norm, we also obtain the boundedness of $M_R$ from $BV(U)$ into
$L^1(U)$ .

\begin{remark}\label{alt} It is possible to define $BV(U)$, where $U$ is open in $\mathbb{R}^d$,
 without knowing a priori
that $|Df|$ is a Radon measure: Write
\begin{equation}\label{defvar}
\int_U |Df| := \sup\left\{\int_U f \operatorname{div} g: g\in C^1_c (U,\mathbb{R}^d), \|g\|_\infty \le 1\right\}.
\end{equation}
Then $f\in BV(U)$ if $f\in L^1(U)$ and $\int_U |Df| < \infty$ (cf., for instance,
Definition 1.3, pg. 4 of \cite{Giu}, or Definition 3.4, pg. 119 and Proposition 3.6,
pg. 120 of \cite{AFP}). Integration by parts immediately yields
that if $f\in C^1(U)$, then
\begin{equation*}
\int_U |Df| = \int_U |\nabla f|dx,
\end{equation*}
(this is Example 1.2 of \cite{Giu}).
With this approach
 one has the following semicontinuity and
approximation results (cf. Theorems 1.9 and 1.17 of \cite{Giu}),
without any reference to Radon measures.
\end{remark}

\begin{theorem}\label{semi} If a sequence of functions $\{f_n\}$ in $BV(U)$ converges in $L^1_{loc}(U)$ to $f$, then $\int_U|Df|\le \liminf_n \int_U |Df_n|$.
\end{theorem}

\begin{theorem}\label{approx} If $f\in BV(U)$, then there exists a sequence of functions $\{f_n\}$ in $BV(U)\cap C^\infty(U)$ such that $\lim_n \int_U |f -f_n| dx = 0$ and
 $\int_U|Df|= \lim_n \int_U |Df_n|$.
\end{theorem}

Note that by passing to a subsequence, we may also assume that $\{f_n\}$ converges
to $f$ almost everywhere.

If one uses the definition of $BV(U)$ given in Remark \ref{alt}, the fact that $Df$ is
a Radon measure is obtained a posteriori via the Riesz Representation Theorem. Then
of course $\int_U |Df| = |Df|(U)$.

\begin{definition} A finite cone $C$ of height $r$, vertex at $0$,
axis $v$, and aperture angle $\alpha$,
is the subset of $B(0,r)$ consisting of all vectors $y$ such that
the angle between $y$ and $v$ is less than or equal to $\alpha /2$.
 A finite cone $C_x$ with vertex at $x$,
is a set of the form $x + C$, where the vertex of $C$ is $0$. Finally, an open set $U$ satisfies a cone condition if there exists a fixed finite cone $C$ such that every
$x\in U$ is the vertex of a cone obtained from $C$ by a rigid motion.
\end{definition}

We shall assume a cone condition in order to have available the following special case of
the  Sobolev embedding theorem (see, for instance, Theorem 4.12, pg. 85 of \cite{AdFo}).
Of course, other type of conditions which also ensure the existence of such an embedding
could be used instead (e.g., $U$ is an extension domain). The next Theorem and its Corollary are well known and included here for the sake of
readability.

\begin{theorem}\label{Sobemb} Let the open set $U\subset \mathbb{R}^d$ satisfy a cone condition. Then
there exists a constant $c>0$, depending only on $U$, such that for all $f\in W^{1,1}(U)$,
$\|f\|_{L^{\frac{d}{d-1}}(U)}\le c \|f\|_{W^{1,1}(U)}$.
\end{theorem}

\begin{corollary}\label{BVemb} Let the open set $U\subset \mathbb{R}^d$ satisfy a cone condition. Then
there exists a constant $c>0$, depending only on $U$, such that for all $f\in BV(U)$,
$\|f\|_{L^{\frac{d}{d-1}}(U)}\le c \|f\|_{BV(U)}$.
\end{corollary}
\begin{proof} Let   $\{f_n\}$ be a sequence of functions in $BV(U)\cap C^\infty(U)$ such that
$f_n\to f$ a.e., $\lim_n \int_U |f -f_n| dx = 0$, and
 $\int_U|Df|= \lim_n \int_U |\nabla f_n|dx$. By Fatou's lemma and Theorem \ref{Sobemb},
$\|f\|_{L^{\frac{d}{d-1}}(U)}\le \liminf_n \|f_n\|_{L^{\frac{d}{d-1}}(U)}
\le \lim_n c \|f_n\|_{W^{1,1}(U)}= c \|f\|_{BV(U)}$.
\end{proof}

The next definition and lemma are valid for an arbitrary set $E\subset \mathbb{R}^k$,
with measure defined by the restriction of the Lebesgue outer measure to the
$\sigma$-algebra of all intersections of Lebesgue
sets with $E$.

\begin{definition}\label{llog} Let $E\subset \mathbb{R}^k$ and $r\ge 1$. A function $g$ belongs to the Banach space $L(\log^+L)^r (E)$ if for
some $t > 0$ we have
\begin{equation}\label{condllogl}
 \int \frac{|g(x)|}{t} \left(\log^+ \frac{|g(x)|}{t}\right)^r dx< \infty.
\end{equation}
 In that case the Luxemburg norm of $g$ is
\begin{equation*}
  \|g\|_{L(\log^+L)^r} := \inf\left\{t > 0: \int \frac{|g(x)|}{t} \left(\log^+ \frac{|g(x)|}{t}\right)^r dx \le 1\right\}.
\end{equation*}
\end{definition}

Note that by monotone convergence the inequality
\begin{equation*}
\int \frac{|g(x)|}{t} \left(\log^+ \frac{|g(x)|}{t}\right)^r dx \le 1
\end{equation*}
holds when $t = \|g\|_{L(\log^+L)^r}$.

We mention that on finite measure spaces, the condition of
Definition \ref{llog} is equivalent to the seemingly stronger
requirement that for all $t> 0$, (\ref{condllogl}) hold.

The next lemma must be well known, but we include it for the reader's convenience.
While stated for all $r\ge 1$, we only need the cases $r=1$ (used  in
Remark \ref{better}), $r=d-1$
(used  in  Theorem \ref{trivialdir}) and $r=d$ (used in Theorem \ref{L1bounds}).

\begin{lemma}\label{logemb} Let $E\subset \mathbb{R}^d$, where $d\ge 2$, and let $r\ge 1$. If $g\in L^{\frac{d}{d-1}}(E)$, then $g\in L(\log^+L)^r(E)$ and   $\|g\|_{L(\log^+L)^r(E)}\le \left(r(d - 1)\right)^{\frac{r(d-1)}d}\|g\|_{L^{\frac{d}{d-1}}(E)}.$
\end{lemma}
\begin{proof} Note that $\log^+ y \le y^\alpha/\alpha$ for
all $y,\alpha>0$, so given $t>0$, if we set  $y = \frac{|g(x)|}{t}$ and $\alpha = \frac{1}{r(d-1)}$,
we get
\begin{equation*}
  \int \frac{|g(x)|}{t} \left(\log^+ \frac{|g(x)|}{t}\right)^r
  dx\le
  \left(r(d - 1)\right)^{r}\left\|\frac{g}{t}\right\|_{L^{\frac{d}{d-1}}(E)}^{\frac{d}{d-1}}.
\end{equation*}
Now let $t_0 < \|g\|_{L(\log^+L)^r}$. Then
$1 < \left(r(d - 1)\right)^{r}\left\|\frac{g}{t_0}\right\|_{L^{\frac{d}{d-1}}(E)}^{\frac{d}{d-1}}$,
from which it follows that $\|g\|_{L(\log^+L)^r(E)}\le \left(r(d - 1)\right)^{\frac{r(d-1)}d}\|g\|_{L^{\frac{d}{d-1}}(E)}.$
\end{proof}

The proof of the next result is similar to that of Theorem \ref{bd}. We indicate the main differences: 1) In Theorem \ref{bd}, since $d=1$, no cone condition appears and we give a fully explicit constant;  2) when $d = 1$, we use the trivial embedding of
$BV(I)$ in $L^\infty$ given in (\ref{eq3}) instead of Corollary \ref{BVemb} and Lemma \ref{logemb}; 3) for $d > 1$, bounds on
the distributional gradient of the corresponding maximal operator are either false or not known.

 \begin{theorem}\label{L1bounds} Let the open set $U\subset \Bbb R^d$ satisfy
a cone condition. For every $R>0$, the local iterated directional maximal operator
$M_R^{d}\circ \dots\circ M_R^{1}$ and the
 local strong maximal operator $M_R^S$  map $BV(U)$ into $L^1(U)$
boundedly. Hence, so do the following operators: The standard local uncentered maximal operator
$M_R$ associated to an arbitrary norm, the local directional maximal operator $M_R^v$, and
$M_R^{i_k}\circ \cdots\circ M_R^{i_1}$, where $1\le k < d$ and $i_1 < \dots <i_k$.
In fact, if  $S_R$ is any of the above
maximal operators, then there exists a constant $c > 0$, which depends only on the open
set $U$, such that for all $f\in BV(U)$,
\begin{equation}\label{conc}
  \|S_Rf\|_{L^1(U)} \le c \left(\|f\|_{BV(U)} + (\log^+R)^d \|f\|_{L^1(U)} \right).
\end{equation}
 \end{theorem}
\begin{proof}
By Corollary \ref{BVemb}, it is enough to
show that
\begin{equation}\label{LplusL}
\|S_Rf\|_{L^1(U)} \le c \left(\|f\|_{L^{d/(d-1)}(U)} + (\log^+R)^d
\|f\|_{L^1(U)} \right).
\end{equation}
Now we can assume that $U=\mathbb{R}^d$. Else, we
extend
 $f$  without changing  the right hand side of  (\ref{LplusL}), by setting $f= 0$  on $\mathbb{R}^d\setminus U$.

The reason we are interested in having
$U=\mathbb{R}^d$ is that later on, we will use the
pointwise equivalence on $\mathbb{R}^d$ of
maximal functions associated to different norms.

By $\eta$ we denote a generic $d$-tuple of integers $(n_1,\ldots,n_d)\in
\mathbb{Z}^d$. For $\eta\in\mathbb{Z}^d$, we define the cubes
$I_{\eta}=[n_1 R,(n_1+1)R)\times\dots\times [n_d R,(n_d +1)R)$ and
$J_{\eta}=[(n_1-1)R,(n_1+2)R)\times\dots\times [(n_d-1)R,(n_d+2)R)$.
Set $f_{\eta}=f|_{J_{\eta}}$.

We want to estimate
\begin{equation*}
\alpha_{\eta}:=  \int_{I_{\eta}}M^d_R\circ\cdots\circ M^1_Rf(x)dx
\end{equation*}
\begin{equation*}
 =
\int_{I_{\eta}}M^d_R\circ\cdots\circ M^1_Rf_{\eta}(x)dx \le
  \int_{I_{\eta}}M^d\circ\cdots\circ M^1f_{\eta}(x)dx.
\end{equation*}
From \cite[\S I. Theorem 1]{Fa}, we get
\begin{equation}\label{weakt}
  \lambda^d(\{M^d\circ\cdots\circ M^1f_{\eta}>4t\})\le C
  \int_{J_{\eta}} \frac{|f_{\eta}(x)|}{t}
  \left(\log^+\frac{|f_{\eta}(x)|}{t}\right)^{d-1} dx,
\end{equation}
where $C$ is a constant that depends only on $d$. Moreover, calling
$A= \|f_{\eta}\|_{L(\log^+L)^{d}}$ and using (\ref{weakt}) we obtain
\begin{equation*}
  \alpha_{\eta} = 4 \int_0^\infty \lambda^d(I_{\eta}\cap \{M^d_R\circ\cdots\circ M^1_Rf_{\eta}(x)>4t\})dt
=4\int_0^{A/R^d} + 4 \int_{A/R^d}^\infty
\end{equation*}
\begin{equation}\label{tres}
\le 4 A + 4C\int_{A/R^d}^\infty \int_{J_{\eta}}
\frac{|f_{\eta}(x)|}{t}
  \left(\log^+\frac{|f_{\eta}(x)|}{t}\right)^{d-1} dx dt = 4A + B.
\end{equation}
Let $\tilde{J}_{\eta}:=J_{\eta}\cap\{|f(x)|>A/R^d\}$. Applying the
Fubini-Tonelli Theorem and the change of variable $y(t)
=\log\frac{|f_{\eta}(x)|}{t}$ we have
\begin{equation*}
 B = 4 C\int_{\tilde{J}_{\eta}}  \int_{A/R^d}^{|f_{\eta}(x)|}
\frac{|f_{\eta}(x)|}{t}
  \left(\log\frac{|f_{\eta}(x)|}{t}\right)^{d-1}dt dx
\end{equation*}
\begin{equation*}
= 4   C\int_{\tilde{J}_{\eta}}
  |f_{\eta}(x)|dx\int_0^{\log^+\frac{|f_{\eta}(x)|R^d}{A}}
  y^{d-1}dy
  \end{equation*}
  \begin{equation*}= \frac{4 C}{d}\int_{\tilde{J}_{\eta}}
  |f_{\eta}(x)|\left(\log\frac{|f_{\eta}(x)|}{A} + d \log R\right)^d  dx
\end{equation*}
  \begin{equation*}\le\frac{4 C 2^d}{d}\int_{{J}_{\eta}}
  |f_{\eta}(x)|\left(\left(\log^+\frac{|f_{\eta}(x)|}{A}\right)^d + d^d \left(\log^+ R\right)^d\right)dx
\end{equation*}
\begin{equation*}
= \frac{4 C 2^d }{d} \left(A\int_{J_{\eta}}
 \frac{ |f_{\eta}(x)|}{A}\left(\log^+ \frac{|f_{\eta}(x)|}{A}\right)^d
  dx + d^d\|f_{\eta}\|_{L^1(J_{\eta})} (\log^+R)^d\right)
  \end{equation*}
  \begin{equation}\label{last}
  \le \frac{4 C 2^d }{d}\left(A + d^d\|f_{\eta}\|_{L^1(J_{\eta})}  (\log^+
  R)^d\right).
\end{equation}
Putting together  (\ref{tres}), (\ref{last}), and Lemma
\ref{logemb}, we get
\begin{equation*}
  \alpha_{\eta}\le C^\prime \left(\|f_{\eta}\|_{L^{d/(d-1)}(J_{\eta})}+ \|f_{\eta}\|_{L^1(J_{\eta})}(\log^+R)^d\right).
\end{equation*}
Next we sum  over  all $d$-tuples $\eta\in  \mathbb{Z}^d$. Since a
point in $\mathbb{R}^d$ cannot be contained in more than $3^d$
different cubes of type $J$, we conclude that for some $c > 0$,
\begin{equation}\label{logd}
  \int_{\mathbb{R}^d} M_R^d\circ \cdots\circ M_R^1f(x) dx \le c \left( \|f\|_{L^{d/(d-1)}(\mathbb{R}^d)}+  \|f\|_{L^1(\mathbb{R}^d)}(\log^+R)^d\right).
  \end{equation}
Since $M_R^Sf(x) \le M_R^d\circ \cdots\circ M_R^1f(x)$ for almost all $x
\in \mathbb{R}^d$, the same inequality holds for $M_R^Sf$. Likewise,
$M_R^S$
 dominates pointwise the maximal operator $M_R$ associated to
the $l^\infty$ norm (i.e., to cubes), so (\ref{conc}) also holds for
$M_R$. Since local maximal operators associated to different norms
are pointwise comparable by the equivalence of all norms in
$\mathbb{R}^d$, inequality (\ref{conc}) holds, perhaps with a
different value of $c$, for the maximal operator $M_R$ defined by
any given norm. Finally, if $1\le k < d$ and $i_1 < i_2 <\dots<i_k$,
we have $M_R^{i_k}\circ \cdots\circ M_R^{i_1}f(x) \le M_R^d\circ
\cdots\circ M_R^1f(x)$ for all $x\in \mathbb{R}^d$, and $M_R^v$
obviously satisfies the same bounds as $M^1_R$, so (\ref{conc})
holds for all the operators under consideration.
\end{proof}

\begin{remark}\label{better} It is possible to obtain  bounds for $M_R$ directly, using essentially
the same proof as in the previous theorem, rather than deriving them from the corresponding
bounds for $M^S_R$.  In fact, a direct approach yields a lower order of growth, $O(\log R)$ instead
of $O((\log R)^d)$. More precisely, replace  in the proof
$L(\log^+L)^d$ by $L(\log^+L)$, and
 inequality (\ref{weakt}) by the following well known refinement (due to N. Wiener, cf. \cite[Theorem 4$^\prime$]{Wi}) of the
weak type inequality:
\begin{equation*}
  \lambda^d(\{Mf>t\})\le \frac{C}{t}\int_{\{|f|>t/2\}}|f(x)|dx
  \qquad \textnormal{for all }t>0.
\end{equation*}
Then argue as before, to get
\begin{equation*}
  \int_U M_R f(x) dx \le c \left( \|f\|_{BV(U)}+  \|f\|_{L^1(U)}\log^+R\right).
  \end{equation*}
  An analogous remark can be made with respect to the operators $M_R^{i_k}\circ \cdots\circ M_R^{i_1}$  and $M_R^v$, obtaining orders of growth $O(\log^k R)$ and $O(\log R)$ respectively.
\end{remark}

\begin{theorem}\label{Strong} Let $d > 1$ and let $U\subset\mathbb{R}^d$ be open.  Given any $R>0$, the following maximal operators are unbounded on $BV(U)$: The local directional
maximal operator $M_R^v$, the local iterated directional maximal operator $M_R^{d}\circ \dots\circ M_R^{1}$, and  the local strong maximal operator $M_R^S$.
 \end{theorem}
\begin{proof} We will show that if $S_R$ denotes any of the maximal operators considered
in the statement of the theorem, then there exists a sequence of characteristic functions
$f_{1/n}$ such that $\lim_{n\to\infty}\|f_{1/n}\|_{BV(U)} = 0$ and $$\lim_{n\to\infty}\frac{|DS_R(f_{1/n})|(U)}{\|f_{1/n}\|_{BV(U)}} =
\infty.$$
  In fact, the same result holds for the corresponding nonlocal maximal operators,
which can be included in the notation $S_R$ by allowing  the possibility $R=\infty$, as
we do in this proof. So we take $0<R\le\infty$. Actually it is enough to consider
$2 < R\le \infty$, since the argument we give below adapts to smaller values for $R$ just by rescaling. Similarly it is enough to consider the case $U=\mathbb{R}^d$. We start with
$M^v_R$. By a rotation we may assume that $v = e_1$. For notational simplicity,
we
will write the proof for the case $d = 2$ only.
Fix $R$. Given $0<\delta< 1$, set $f_\delta(x):=\chi_{[0,\delta]^2}(x)$.
Then
\begin{equation*}
\|f_\delta\|_1=\delta^2
\end{equation*}
and, since  $|Df_\delta|(\mathbb{R}^2)$ is just the perimeter of the square $[0,\delta]^2$
(cf., for instance, Exercise 3.10 pg. 209 of \cite{AFP}),
\begin{equation*}
|Df_\delta|(\mathbb{R}^2)=4\delta.
\end{equation*}
Thus
\begin{equation}\label{une}
\|f_\delta\|_{BV(\mathbb{R}^2)} = O(\delta) \textnormal{ when }
\delta \rightarrow 0.
\end{equation}
Next, let $\delta\le x\le 1$, and $0\le y\le \delta$. It is then easy to check that
\begin{equation*}
M_R^1(f_\delta)(x,y)=\frac{\delta}{x}.
\end{equation*}
Given $\delta \le t < 1$, the level sets $E_t :=\{M_R^1 (f_\delta ) > t\}$ are rectangles,
with perimeter
$$|D\chi_{E_t}|(\mathbb{R}^2) \ge 2\delta + \frac{2\delta}{t}.$$
 By the coarea formula for BV functions
(cf. Theorem 3.40, pg. 145 of \cite{AFP}), we have
\begin{equation}\label{coar}
|DM^1_R f_\delta|(\mathbb{R}^2)=\int_{-\infty}^{\infty}
|D\chi_{E_t}|(\mathbb{R}^2) dt \ge \int_{\delta}^{1}
|D\chi_{E_t}|(\mathbb{R}^2) dt \ge 2\delta \int_{\delta}^{1}\left( 1
+ \frac{1}{t}\right) dt = \Theta\left(\delta\log
\frac1{\delta}\right).
\end{equation}
where $\Theta$ stands for the exact order of growth.
 From (\ref{une}) and (\ref{coar}) we obtain
\begin{equation}\label{unbound}
\frac{|DM_R^1(f_\delta)|(\mathbb{R}^2)}{\|f_\delta\|_{BV((\mathbb{R}^2)}}\rightarrow
\infty \textnormal{ when }\delta\rightarrow 0,
\end{equation}
as was to be proven.

Note next that on $[0,1]\times [0,\delta]$ the three maximal functions $M^1_Rf_\delta$,
 $M_R^2\circ M_R^1 f_\delta$ and $M_R^S f_\delta$ take the same values, from which
 it easily follows that for $\delta \le t < 1$,
$$|D\chi_{\{M_R^2\circ M_R^1 (f_\delta ) > t\}}|(\mathbb{R}^2) \ge 2\delta + \frac{2\delta}{t}$$
and
$$|D\chi_{\{M_R^S (f_\delta ) > t\}}|(\mathbb{R}^2) \ge 2\delta + \frac{2\delta}{t}.$$
Thus,  the analogous statement to (\ref{unbound})  holds for  $M_R^2\circ M_R^1 f_\delta$ and $M_R^S f_\delta$ also.
\end{proof}

A standard mollification argument shows that the preceding maximal operators are
not bounded on $W^{1,1}(U)$ either.

\section{Converses and a one dimensional characterization.}

 Recall that   $f^+$ and $ f^-$  denote respectively the positive and negative parts of $f$.
Now, for  any open set $U\subset \Bbb R^d$,  $f\in BV(U)$  if and
only if  both $f^+\in BV(U)$ and $f^-\in BV(U)$. This can be seen as
follows: If $f\in BV(U)$, it is immediate from the definition
\ref{defvar} contained in Remark \ref{alt} that $\int_U |Df| \ge
\int_U |D(f^+)|$ and  $\int_U |Df| \ge \int_U |D(f^-)|$, so $f^+,
f^-\in BV(U)$.
 On the other hand, if both $f^+, f^-\in
BV(U)$, then there are  sequences $\{g_n\}$ and $\{h_n\}$ of
$C^\infty$ functions that approximate  $f^+$ and $f^-$ respectively,
in the sense of Theorem \ref{approx}. Since $g_n - h_n\to f$ in
$L^1(U)$, by semicontinuity $|Df|(U) \le \liminf_n \int_U |\nabla
(g_n - h_n)| dx \le \lim_n \int_U |\nabla g_n| dx +\lim_n \int_U
|\nabla h_n| dx =  |D(f^+)|(U) + |D(f^-)(U)|.$ Hence $f\in BV(U)$.

\begin{theorem}\label{trivialdir}
    Let  $U\subset \Bbb R^d$ be an open set
and let  $f:U\to \Bbb R $ be locally integrable. Suppose that there
exists a sequence $\{a_n\}_1^\infty$ with $\lim_n a_n = 0$ and a
constant $c$ such that for all $n$,
   $M_{a_n}f^+\in W^{1,1}(U)$, $M_{a_n}f^-\in W^{1,1}(U)$, $\| M_{a_n}f^+\|_{W^{1,1}(U)} \le c$, and $\| M_{a_n}f^-\|_{W^{1,1}(U)} \le c$.
Then $f\in BV(U)$. The same happens if instead of $M_R$ we consider
either the local directional maximal operator, or,  under the
additional hypothesis that $U$ satisfies a cone condition, the local
strong
 maximal operator.
\end{theorem}

\begin{proof} Consider first $f^+$. By the Lebesgue Theorem on differentiation of integrals we have that
$\lim_n M_{a_n} f^+ = f^+$ a.e., so by dominated convergence, $M_{a_n} f^+ \to  f^+$ in
$L^1(U)$, and by Theorem \ref{semi}, $\int_U|Df^+| \le \liminf_n
\int_U |DM_{a_n} f^+| \le c <\infty$. Repeating the argument for $f^-$ we get
$|Df|(U) \le |Df^+|(U) + |Df^-|(U) <\infty$.
The result for the local strong  maximal operator follows from
 the well known Theorem of Jessen, Marcinkiewicz  and Zygmund (\cite{JMZ}) stating that
basis of rectangles (with sides parallel to the axes)
 differentiates  $L(\log ^+ L)^{d-1}_{loc}(U)$, and hence $BV(U)$ (cf. Corollary \ref{BVemb} and Lemma \ref{logemb}; for the first embedding we
use the cone condition). Finally, the weak type $(1,1)$ boundedness
of $M_T^v$ (which is obtained from the one dimensional result and
the Fubini-Tonelli Theorem) also entails, by the standard argument,
the corresponding differentiation of integrals result, so $\lim_n
M^v_{a_n} f^+ = f^+$ and $\lim_n M^v_{a_n} f^- = f^-$.
\end{proof}

For intervals $I\subset \mathbb{R}$ we have the following characterization.

\begin{theorem}\label{charact}
    Let   $f:I\to \Bbb R $ be locally integrable. Then the following are equivalent:

a) $f\in BV(I )$.

b) $M_Rf^+\in W^{1,1}(I)$, $M_Rf^-\in W^{1,1}(I)$,
$\|M_Rf^+\|_{W^{1,1}(I)}\le 3(1+2\log^+(R)) \|f^+\|_{L^1(I)} + 4 |Df^+|(I),$ and $\|M_Rf^-\|_{W^{1,1}(I)}\le 3(1+2\log^+(R)) \|f^-\|_{L^1(I)} + 4 |Df^-|(I).$

c) There exists a sequence $\{a_n\}_1^\infty$ with $\lim_n a_n = 0$  and a constant
 $c = c(f, \{a_n\}_1^\infty)$ such that for all $n$,
   $M_{a_n}f^+\in W^{1,1}(I)$, $M_{a_n}f^-\in W^{1,1}(I)$, $\| M_{a_n}f^+\|_{W^{1,1}(I)} \le c$, and $\| M_{a_n}f^-\|_{W^{1,1}(I)} \le c$.

d) There exists an $R>0$ and a constant $c = c(f, R)$ such that for all
$T\in (0,R]$,   $M_Tf^+\in W^{1,1}(I)$, $M_Tf^-\in W^{1,1}(I)$, $\| M_Tf^+\|_{W^{1,1}(I)} \le c$, and $\| M_Tf^-\|_{W^{1,1}(I)} \le c$.

e) For every $R>0$ there exists a constant $c = c(f,R)$ such that for all
$T\in (0,R]$,   $M_Tf^+\in W^{1,1}(I)$, $M_Tf^-\in W^{1,1}(I)$, $\| M_Tf^+\|_{W^{1,1}(I)} \le c$, and $\| M_Tf^-\|_{W^{1,1}(I)} \le c$.

 If   $f:I\to \Bbb R $ is absolutely continuous, then
 a') $f\in W^{1,1}(I)$ is equivalent to b), c), d) and e).
 \end{theorem}

\begin{proof} The implications b) $\to $ e),  e) $\to $ d) and d) $\to $ c) are obvious,
and a) $\to$ b) is the content of  Theorem \ref{bd}. Without
loss of generality we may take $I$ to be open, so c) $\to $ a) is a special case
of Theorem \ref{trivialdir}. Finally, the last claim follows from the fact that
 $f\in W^{1,1}(I)$  if and only if $f$ is absolutely continuous and
 $f\in BV(I)$.
\end{proof}

 Let   $f:I\to \Bbb R $ be locally integrable.
By Theorem \ref{bd}, if  $|f|\in BV(I )$ then for every $R >0$,
 $M_Rf\in W^{1,1}(I)$ boundedly, with bound depending on $R$. Thus it is natural to ask whether
 the latter condition alone suffices to ensure that  $|f|\in BV(I )$. In other words, we
 are asking whether the uniform bound condition appearing in parts c),  d) and e) of
 Theorem \ref{charact} is really
needed. The following example shows that the answer is positive.

\begin{example}\label{counter1d}{\em There exists a non-negative function $f\in L^1(\Bbb R)\setminus BV(\Bbb R)$ such  that for all
$R>0$,  $M_Rf\in W^{1,1}(\Bbb R )$.}
\begin{proof} Let $A$ be the closed set $[-1000,0]\cup\left(\cup_{n=0}^\infty [2^{-n}, 2^{-n} + 2^{-n-1}]\right)$,
and let $f$ be the upper semicontinuous function $\chi_A$. Fix $R>0$. Clearly $M_R f \ge f$ everywhere,
so by Lemma 3.4 of \cite{AlPe}, $M_Rf$ is a continuous function. Also, $M_R f|_{\mathbb{R}\setminus (0,2^{-n})}$ is Lipschitz, with $\operatorname{Lip }(M_Rf)\le \max\{R^{-1}, 2^{n+1}\}$, by Lemma 3.8 of \cite{AlPe}. Hence, if $E\subset \Bbb R$ has measure zero, so does $M_Rf(E)$, being a countable union of sets of measure zero. Next we show that $|DM_R f|(\Bbb R) < \infty$.  Let $n\ge 1$ . On intervals of the
form $( 2^{-n} + 2^{-n-1},  2^{-n+1})$, if $R > 2^{-n-2}$ then $M_Rf > f$, so by Lemma 3.6 of \cite{AlPe}
there exists an
$x_n\in ( 2^{-n} + 2^{-n-1},  2^{-n+1})$ such that
  $M_Rf$ is decreasing on $( 2^{-n} + 2^{-n-1},x_n)$ and increasing on $(x_n, 2^{-n+1})$.
Taking this fact into account, it is easy to see
 that $V(M_R f, \mathbb{R})$ is decreasing in $R$, so we may suppose $R \in (0,1)$.
Select $N\in \Bbb N$ such that $2^{-N+1} < R$. Then for $n > N$,
\begin{equation*}
    V(M_Rf,( 2^{-n} + 2^{-n-1},  2^{-n+1}))
= 2\left(1-M_Rf(x_n)\right)
\le
2\left(1-\frac{R- 2^{-n+1}}{R}\right)\le \frac{ 2^{-n+2}}{R}.
\end{equation*}
Hence $|DM_R f|(\Bbb R ) \le 2 + 2(N+1) <\infty$.
Since $M_Rf$  is continuous, of bounded variation, and maps measure
zero sets into measure zero sets, by the Banach Zarecki Theorem it
 is absolutely continuous, so  $M_Rf\in W^{1,1}(\Bbb R )$.
\end{proof}
Of course, using $\Bbb R$ above is not necessary, the example can be easily adapted to any
other interval $I$.
\end{example}

\end{document}